\documentclass[12pt]{amsart}
\usepackage[top=1in,left=1in,bottom=1in,right=1in]{geometry}
\usepackage{amsmath}
\usepackage{amssymb}
\usepackage{amsthm}
\usepackage{tikz}
\usetikzlibrary{arrows}
\usepackage{graphicx}
\usepackage{url}
\usepackage{colordvi}
\usepackage{color}
\usepackage{verbatim}

\newcommand{\ZZ}[1]{\mathbb{Z}/#1\mathbb{Z}}

\newtheorem{theorem}{Theorem}[section]
\newtheorem{lemma}[theorem]{Lemma}
\newtheorem{conjecture}[theorem]{Conjecture}

\newtheorem{problem}[theorem]{Problem}
\newtheorem{proposition}[theorem]{Proposition}
\newtheorem{corollary}[theorem]{Corollary}
\theoremstyle{definition}
\newtheorem{definition}[theorem]{Definition}
\theoremstyle{remark}
\newtheorem{remark}[theorem]{Remark}
\newtheorem{example}[theorem]{Example}

\begin{document}
\title{Chromatic Bounds on Orbital Chromatic Roots}

\author{Dae Hyun Kim}
\address{Department of Mathematics.  California Institute of Technology, Pasadena, CA, 91125}
\email{dkim3@caltech.edu}

\author{Alexander H. Mun}
\address{Department of Mathematics.  California Institute of Technology, Pasadena, CA, 91125}
\email{amun@caltech.edu}

\author{Mohamed Omar}
\address{Department of Mathematics.  Harvey Mudd College, Claremont, CA 91711}
\email{omar@g.hmc.edu}



\begin{abstract}
Given a group $G$ of automorphisms of a graph $\Gamma$, the orbital chromatic polynomial $OP_{\Gamma,G}(x)$ is the polynomial whose value at a positive integer $k$ is the number of orbits of $G$ on proper $k$-colorings of $\Gamma.$  Cameron and Kayibi introduced this polynomial as a means of understanding roots of chromatic polynomials.  In this light, they posed a problem asking whether the real roots of the orbital chromatic polynomial of any graph are bounded above by the largest real root of its chromatic polynomial. We resolve this problem in a resounding negative by not only constructing a counterexample, but by providing a process for generating families of counterexamples. We additionally begin the program of finding classes of graphs where the answer to this problem is true; in particular establishing its veracity for many outerplanar graphs.
\end{abstract}

\thanks{The first and second authors were supported by Summer Undergraduate Research Fellowships at the California Institute of Technology.}
\date{\today}
\maketitle

\section{Introduction}

The chromatic polynomial of a graph $\Gamma$, denoted $P_{\Gamma}(x)$, is the function whose value at any positive integer $k$ is the number of proper $k$-colorings of $\Gamma$.  That $P_{\Gamma}(x)$ is indeed a polynomial comes directly from the classical fact that it satisfies a deletion-contraction relation.  Chromatic polynomials were introduced by Birkhoff \cite{Birkhoff} in 1912, with the particular intent of algebraically resolving what was then the 4-Color Conjecture; indeed this amounts to establishing that any planar graph $\Gamma$ satisfies $P_{\Gamma}(4)>0$.  This perspective motivated the algebraic study of the roots of chromatic polynomials in general.

Though Birkhoff did not manage to prove the 4-Color Conjecture algebraically, Birkhoff and Lewis \cite{BirkhoffLewis} did prove that for planar graphs $\Gamma$, $P_{\Gamma}(x)>0$ for $x \in [5,\infty)$. They additionally conjectured the still open problem that when $\Gamma$ is planar, $P_{\Gamma}(x)>0$ for $x \in [4,\infty)$.  The interval in this conjecture can not be extended, as proven by Royle \cite{Royle}, where it is shown that real chromatic roots can come arbitrarily close to $4$.  Outside of the context of planarity, Sokal \cite{Sokal} proved that the complex roots of chromatic polynomials are dense in the complex plane.  However, in contrast to this, Jackson (see \cite{Jackson}) proved proved $(1,\frac{32}{27}]$ is a zero-free real interval for chromatic roots, and Thomassen (see \cite{Thomassen}) surprisingly proved that real chromatic roots are dense in $[\frac{32}{27},\infty)$.

Though investigating algebraic properties of chromatic polynomials has been fruitful, they distinguish between colorings that can be obtained from one another through an automorphism of a graph.  This motivates the definition of the \emph{orbital chromatic polynomial} (as introduced in \cite{Cameron}), which does not distinguish between two such colorings.  In particular, given a graph $\Gamma$ and a group of automorphisms $G$ of $\Gamma$, the orbital chromatic polynomial $OP_{\Gamma,G}(x)$ is the function whose value at a positive integer $k$ is the number of $G$-orbits of proper $k$-colorings of $\Gamma$.  

 That $OP_{\Gamma,G}(x)$ is indeed a polynomial in $x$ can be seen as follows.  For any $g \in G$, define $P_{\Gamma/g}(k)$ to be the number of proper $k$-colorings of $\Gamma$ fixed by $g$.  By the Orbit-Stabilizer Theorem,
\begin{equation}\label{OP}
OP_{\Gamma,G}(k) = \frac{1}{|G|} \sum_{g \in G} P_{\Gamma/g}(k).
\end{equation}
The values $P_{\Gamma/g}(k)$ are in fact evaluations of chromatic polynomials themselves.  To see this, construct the graph which we conveniently name $\Gamma/g$, whose vertices are the orbits of the action of $\langle g \rangle$ on $\Gamma$, with two orbits $O_1,O_2$ adjacent if there are vertices $v_1 \in O_1$ and $v_2 \in O_2$ such that $v_1$ is adjacent to $v_2$ in $\Gamma$.  Since a coloring is fixed by $g$ if and only if it is constant on orbits of the action of $g$, any proper coloring of $\Gamma$ fixed by $g$ induces a proper coloring on $\Gamma/g$ and vice-versa.  As a result, for any positive integer $k$, $P_{\Gamma/g}(k)$ is the number of proper $k$-colorings of $\Gamma/g$.   But this is a polynomial in $k$, so by Equation~(\ref{OP}), $OP_{\Gamma,G}(k)$ is a polynomial in $k$ for any positive integer $k$.  This then implies $OP_{\Gamma,G}(x)$ is a polynomial in $x$, and indeed $OP_{\Gamma,G}(x) = \frac{1}{|G|} \sum_{g \in G} P_{\Gamma/g}(x)$.  Note that this also provides an algorithm for computing $OP_{\Gamma,G}(x)$.

Analogous to studying algebraic properties of chromatic polynomials, of particular interest is understanding algebraic properties of orbital chromatic polynomials, and identifying how phenomena that hold for chromatic polynomials transfer to orbital chromatic polynomials.  This topic is the central focus of \cite{Cameron}.  One of the main results there, in contrast to Jackson's and Thomassen's results on zero-free intervals and the density of real chromatic roots in $[\frac{32}{27},\infty)$, is that real orbital chromatic roots are dense in $\mathbb{R}$.  This proof was constructive, and in all graphs constructed in the proof it was observed that the real roots of $OP_{\Gamma,G}(x)$ were always bounded above by the largest real root of $P_{\Gamma}(x)$.  This led to the following natural problem originally posed in \cite{Cameron}.

\begin{problem}[\cite{Cameron}, Problem 2]\label{prob:Cameron}
Is it true that the real roots of $OP_{\Gamma,G}(x)$ are bounded above by the largest real root of $P_{\Gamma}(x)$ for any graph $\Gamma$ and any subgroup $G$ of $Aut(\Gamma)$?
\end{problem}

We show that the answer is a resounding no.  Indeed, we not only find examples where the property described in Problem~\ref{prob:Cameron} is false, but we additionally provide a means for generating families of graphs for which the property is false, as encapsulated in the following theorem: 

\begin{theorem}\label{counterexample}
Let $\Gamma$ be a graph and $G$ be a group of automorphisms of $\Gamma$.  Suppose the following hold:

\begin{enumerate}
\item There is some $g \in G$ for which $\Gamma/g$ contains fewer vertices than any of the graphs $\{\Gamma/h  :  h \in G, h \neq g\}$.  

\item For the $g$ in part (1), there is some $x_0 \notin \mathbb{Z}$ greater than the largest real root of $P_{\Gamma}(x)$ such that $P_{\Gamma/g}(x_0)<0$.

\end{enumerate}
Then one can construct, from $\Gamma$, a graph $\Gamma'$ and a group of automorphisms of $\Gamma'$, say $G'$, such that $OP_{\Gamma',G'}(x)$ has a real root larger than any real root of $P_{\Gamma'}(x)$.
\end{theorem}
We refer the reader to Example~\ref{ex2} for an illustration of Theorem~\ref{counterexample} and a subsequent example establishing the first negative example to question posed in Problem~\ref{prob:Cameron}.

It now remains to characterize the pairs $(\Gamma,G)$ where $\Gamma$ is a graph and $G$ is a group of automorphisms of $\Gamma$ for which the property described in Problem~\ref{prob:Cameron} holds.  We begin this program by focusing on planar graphs, as these graphs served as the motivation for studying chromatic and orbital chromatic polynomials in the first place.  In this light, we uncover a family of planar graphs for which the property indeed holds:

\begin{theorem}\label{outerplanar}
If $\Gamma$ is an outerplanar graph that contain at least one odd cycle, then the real roots of $OP_{\Gamma,G}(x)$ are bounded above by the largest real root of $P_{\Gamma}(x)$ for any group $G$ of automorphisms of $\Gamma$.
\end{theorem}

\subsection*{Outline}

The organization of the paper is as follows:  In Section 2, we focus on proving Theorem~\ref{counterexample} in order to supply the machinery for finding negative examples to Problem~\ref{prob:Cameron}.  In Section 3, we explore when the property described in Problem~\ref{prob:Cameron} holds, ultimately leading to the proof of Theorem~\ref{outerplanar}. We conclude with some conjectures in Section 4.

\section{Failure: Bounding Orbital Chromatic Roots}

This section is dedicated to proving Theorem~\ref{counterexample} and subsequently constructing negative examples to Problem~\ref{prob:Cameron}.  In order to do this, we need to define some auxiliary graphs.

\begin{definition}\label{kns}
Given positive integers $n,s$, define the graph $K_n$ to be the complete graph on $n$ vertices, and $N_s$ to be the graph consisting of $s$ isolated vertices.  Define the graph $H_{n,s}$, the \emph{join} of $K_n$ and $N_s$, to be the graph obtained by taking the union of $K_n$ and $N_s$, and adding an edge between every vertex in $K_n$ and every vertex in $N_s$.
\end{definition}

The graphs $H_{n,s}$ are central to constructing negative to Problem~\ref{prob:Cameron}.  In particular, the following construction will be crucial:

\begin{definition}\label{suspend}
Let $\Gamma$ be a graph with vertex set $\{v_1,v_2,\ldots,v_k\}$, and $n,s$ be positive integers.  Let $H_{n,s}^{(1)},H_{n,s}^{(2)},\ldots,H_{n,s}^{(k)}$ be $k$ copies of the graph $H_{n,s}$ and choose vertices $u_i \in V(H_{n,s}^{(i)})$ so that there is an isomorphism from $H_{n,s}^{(i)}$ to $H_{n,s}^{(j)}$ sending $u_i$ to $u_j$.  We construct the graph $\Gamma^{(n,s)}$ by starting with $\Gamma$, and appending the $k$ copies of $H_{n,s}$ to $\Gamma$ by identifying the vertices $u_i$ and $v_i$ for $i \in \{1,2,\ldots,k\}$.
\end{definition}
The following proposition is immediate.
\begin{proposition}\label{prop}
For any graph $\Gamma$, and positive integers $n,s$,
\[
P_{\Gamma^{(n,s)}}(x) = \left( (x-1) \cdots (x-n+1)(x-n)^s \right)^{|V(\Gamma)|} \cdot P_{\Gamma}(x) = \left( \frac{P_{H_{n,s}}(x)}{x} \right)^{|V(\Gamma)|} \cdot P_{\Gamma}(x).
\]
\end{proposition}

With these constructions, we can now prove Theorem~\ref{counterexample}.

\proof (of Theorem~\ref{counterexample})
Let $n,s$ be positive integers (arbitrary for now).  Construct the graph $\Gamma^{(n,s)}$ and let $G^{(n,s)}$ be the group induced by $G$ that permutes vertices of the subgraph $\Gamma$ of $\Gamma^{(n,s)}$ just as $G$ does, so that if $g \in G$ sends $v_i$ to $v_j$, then $H_{n,s}^{(i)}$ gets sent to $H_{n,s}^{(j)}$ via the isomorphism sending $u_i$ to $u_j$.  There is a natural bijection between elements in $G$ and elements in $G^{(n,s)}$, so for any $h$ in $G$, we denote by $h^{(n,s)}$ its corresponding element in $G^{(n,s)}$.  

Observe that for any $h \in G$, $\Gamma^{(n,s)}/h^{(n,s)} = (\Gamma/h)^{(n,s)}$, and so by Proposition~\ref{prop}
\[
P_{\Gamma^{(n,s)}/h^{(n,s)}}(x) = P_{(\Gamma/h)^{(n,s)}}(x)= \left( \frac{P_{H_{n,s}}(x)}{x} \right)^{|V(\Gamma/h)|} \cdot P_{\Gamma/h}(x),
\]
and hence $OP_{\Gamma^{(n,s)},G^{(n,s)}}(x)$ is
\[
 \frac{1}{|G|} \left( \frac{P_{H_{n,s}}(x)}{x} \right)^{|V(\Gamma/g)|} \left( P_{\Gamma/g}(x) + \sum_{h \in G, h \neq g} \left( \frac{P_{H_{n,s}}(x)}{x} \right)^{|V(\Gamma/h)| - |V(\Gamma/g)|} P_{\Gamma/h}(x) \right).
\]

We can now choose appropriate values of $n$ and $s$ to control the roots of $OP_{\Gamma^{(n,s)},G^{(n,s)}}(x)$.  First, recall our assumption that there is some $x_0 \notin \mathbb{Z}$ for which $P_{\Gamma/g}(x_0)<0$, and that this $x_0$ is larger than any real root of $P_{\Gamma}(x)$.  Choose $n = \lfloor x_0 \rfloor$.  Since $x_0 \in (n,n+1)$, as we increase $s$ the quantity $\frac{P_{H_{n,s}}(x_0)}{x_0}$ will be positive and approach $0$.  Together with the fact that $P_{\Gamma/g}(x_0)<0$, and that $|V(\Gamma/h)| - |V(\Gamma/g)| > 0$ for all $h \neq g$, this implies we can choose a sufficiently large value of $s$ for which $OP_{\Gamma^{(n,s)},G^{(n,s)}}(x_0) < 0$.  But $\lim_{x \to \infty} OP_{\Gamma^{(n,s)},G^{(n,s)}}(x) = \infty$ so by the Intermediate Value Theorem $OP_{\Gamma^{(n,s)},G^{(n,s)}}(x)$ has a root larger than $x_0$.  

We now construct $\Gamma'$ and $G'$, letting $\Gamma' = \Gamma^{(n,s)}$ and $G' = G^{(n,s)}$ for our particular choices of $n$ and $s$ above.  Then $OP_{\Gamma',G'}(x)$ has a real root larger than $x_0$, whereas \[P_{\Gamma'}(x) = P_{\Gamma}(x) \cdot \left( \frac{P_{H_{n,s}}(x)}{x} \right)^{|V(\Gamma)|},\] whose maximum real root does not exceed $x_0$. \qed

\begin{example}\label{ex2}
Let $\Gamma$ be a $6$-cycle with vertices labeled $\{1,2,3,4,5,6\}$, where the neighbors of $i$ are $i-1$ and $i+1$ (taken mod $6$) for each $i \in V(\Gamma)$.  Let $g$ be the automorphism that sends $i$ to $i+3$ (taken mod $6$), and $G$ be the two-element group consisting of $g$ and the identity $e$.  First, note that $\Gamma/e = \Gamma$ and $\Gamma/g$ is a $3$-cycle, so $|V(\Gamma/g)|<|V(\Gamma/e)|$.  Moreover, observe that
\[
P_{\Gamma}(x) = (x-1)^6+(x-1), \ \ P_{\Gamma/g}(x) = x(x-1)(x-2),
\]
so $x_0=\frac{3}{2}$ is greater than any real root of $P_{\Gamma}(x)$ and $P_{\Gamma/g}(x_0)<0$.  According to the proof of Theorem~\ref{counterexample}, this means we should choose $n=\lfloor x_0 \rfloor = 1$.  In this case, we have that for any $s$,
\[
OP_{\Gamma^{(1,s)},G^{(1,s)}}(x) = \frac{1}{2} (x-1)^{3s} \left( x(x-1)(x-2) + (x-1)^{3s} \left( (x-1)^6 + (x-1) \right) \right),
\]
and hence
\[
OP_{\Gamma^{(1,s)},G^{(1,s)}}\left(\frac{3}{2}\right) = \left(\frac{1}{2}\right)^{3s+1} \cdot \left( -\frac{3}{8} +  \frac{33}{64} \left(\frac{1}{2}\right)^{3s} \right).
\]
Letting $s=1$ we see $OP_{\Gamma^{(1,1)},G^{(1,1)}}\left(\frac{3}{2}\right)<0$ and hence $OP_{\Gamma^{(1,1)},G^{(1,1)}}(x)$ has a real root greater than $\frac{3}{2}$, which is greater than the real roots of $P_{\Gamma^{(1,1)}}(x)$.  See Figure 1 for an illustration of the pertinent graphs in question.
\end{example}

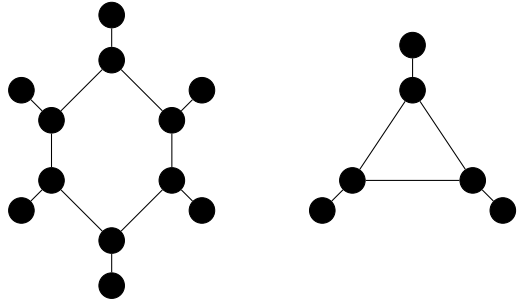
\begin{figure}
\begin{tikzpicture}
  [scale=.4,auto=left,every node/.style={circle, draw, fill=black,
                        inner sep=0pt, minimum width=4pt}]
  \node (n1) at (6,6) {6};
  \node (n2) at (8,4) {5};
  \node (n3) at (8,2) {4};
  \node (n4) at (6,0) {3};
  \node (n5) at (4,2) {2};
  \node (n6) at (4,4) {1};
	\node (n7) at (6,7.5)  {1};
  \node (n8) at (9,5)    {2};
  \node (n9) at (9,1)    {3};
  \node (n10) at (6,-1.5)  {9};
  \node (n11) at (3,1)   {8};
  \node (n12) at (3,5) {7};
	
	\node (n13) at (16,5) {6};
  \node (n23) at (18,2) {5};
  \node (n33) at (14,2) {4};
  \node (n43) at (16,6.5) {3};
  \node (n53) at (19,1) {2};
  \node (n63) at (13,1) {1};

  \foreach \from/\to in {n6/n5,n5/n4,n4/n3,n3/n2,n2/n1,n1/n6,n1/n7,n2/n8,n3/n9,n4/n10,n5/n11,n6/n12,n13/n23,n13/n33,n23/n33,n13/n43,n23/n53,n33/n63}
    \draw (\from) -- (\to);
\end{tikzpicture}

\caption{$\Gamma^{(1,1)}$ (on the left) where $\Gamma$ is a $6$-cycle.  $\Gamma^{(1,1)}/g^{(1,1)}$ (on the right) where $g$ is the $180$-degree rotational symmetry of $\Gamma$.  $\Gamma^{(1,1)}$ is the first of many negative examples to Problem~\ref{prob:Cameron}.}
\end{figure}

\begin{remark}
Note that for any $s \geq 1$, $OP_{\Gamma^{(1,s)},G^{(1,s)}}\left(\frac{3}{2}\right)<0$, so we get a family of counterexamples arising from the graphs $\Gamma^{(1,s)}$ and their group of automorphisms $G^{(1,s)}$ for every positive integer $s$.
\end{remark}

\begin{remark}
Notice that our choice of appending copies of $H_{n,s}$ to each vertex of $\Gamma$ was employed to ensure $\lim_{s \to \infty} \frac{P_{H_{n,s}}(x_0)}{x_0} = 0$, as this was the crux of the argument for finding a graph $\Gamma'$ with a group of automorphisms $G'$ giving rise to a negative example to the property described in Problem~\ref{prob:Cameron}.  We could have easily replaced $H_{n,s}$ with any family of graphs $\{H'_{n,s}\}$ parameterized by natural numbers $n,s$ for which $\lim_{s \to \infty} \frac{P_{H'_{n,s}}(x_0)}{x_0} = 0$ with $n =\lfloor x_0 \rfloor$.  Any such graph family would generate entirely new classes of counterexamples.
\end{remark}

\section{Success: Bounding Orbital Chromatic Roots for Graph Families}

In this section, we start the program of determining which graphs have the property described in Problem~\ref{prob:Cameron}.  We begin by showing this property holds for paths and cycles.  This sets the stage for working our way up to proving Theorem~\ref{outerplanar}, showing that the property described in Problem~\ref{prob:Cameron} holds for outerplanar graphs that have at least one odd cycle.

\subsection{Preliminaries}

We begin by introducing technical preliminaries.  The first of these is the Reduction Lemma, which allows us to conclude that the real roots of $OP_{\Gamma,G}(x)$ are bounded above by the largest real root of $P_{\Gamma}(x)$ if we have certain bounds on the real roots of the chromatic polynomials $\{P_{\Gamma/g}(x) : g \in G\}$.  Though Lemma~\ref{lem:reduction} states the Reduction Lemma in fully generality, we will use Corollary~\ref{cor:reduction} more often in practice.  

\begin{lemma}\label{lem:reduction}(Reduction Lemma)
Let $\Gamma$ be a graph without loops, and let $G$ be a group of automorphisms of $\Gamma$.  Let $\mathcal{G}$ be a partition of $G$, and define $\tilde{\mathcal{G}}$ to consist of those sets $X$ in the partition $\mathcal{G}$ for which $\sum_{g \in X} P_{\Gamma/g}(x) \neq 0$.  If the real roots of the polynomials $\{\sum_{g \in X} P_{\Gamma/g}(x) : X \in \tilde{\mathcal{G}}\}$ are bounded above by the largest real root of $P_{\Gamma}(x)$, then the real roots of $OP_{\Gamma,G}(x)$ are bounded above by the largest real root of $P_{\Gamma}(x).$
\end{lemma}

\begin{corollary}\label{cor:reduction}
Let $\Gamma$ be a graph without loops, and $G$ be any group of automorphisms of $\Gamma$.  Suppose that for all $g$ for which $\Gamma/g$ has no loops, the real roots of $P_{\Gamma/g}(x)$ are bounded above by the largest real root of $P_{\Gamma}(x)$.  Then the real roots of $OP_{\Gamma,G}(x)$ are bounded above by the largest real root of $P_{\Gamma}(x)$.
\end{corollary}
\proof(of Corollary~\ref{cor:reduction}) Apply Lemma~\ref{lem:reduction} where $G$ is partitioned into its individual elements.  The elements $g \in G$ for which $P_{\Gamma/g}(x)$ is non-zero are precisely the ones for which $\Gamma/g$ has no loops. \qed

\proof (of Lemma~\ref{lem:reduction})
For simplicity, define $P_X(x) := \sum_{g \in X} P_{\Gamma/g}(x)$, and let $r$ be the largest real root of $P_{\Gamma}(x)$.  Suppose there is some $r'>r$ that is a root of $OP_{\Gamma,G}(x)$.  Then
\[
0=OP_{\Gamma,G}(r') = \frac{1}{|G|} \sum_{X \in \tilde{\mathcal{G}}} P_X(r').
\] This implies $P_X(r')<0$ for some $X \in \tilde{\mathcal{G}}$.  Since $\lim_{x \to \infty} P_X(x) = \infty$ and $P_X(x)$ is continuous, the Intermediate Value Theorem implies $P_X(x)$ has a real root larger than $r'$.  This contradicts that the largest real root of $P_X(x)$ is at most $r$.  \qed

Another construction that we use throughout the paper is adding path ears to graphs.  In particular, given a graph $\Gamma$, a pair of adjacent vertices $u,v \in V(\Gamma)$, and a positive integer $n$, we denote by $\Gamma_{u,v}(n)$ the graph obtained from $\Gamma$ by adding a path from $u$ to $v$ with $n$ interior vertices, none of which are in $\Gamma$.  Using deletion and contraction, one can inductively prove

\begin{lemma}\label{lem:cyclereduction}
\[
P_{\Gamma_{u,v}(n)}(x) =  \left( \sum_{i=0}^n (-1)^{n-i} (x-1)^i \right) \cdot P_{\Gamma}(x) = (-1)^n \left( \frac{1-(1-x)^{n+1}}{x}\right) \cdot P_{\Gamma}(x).
\]
\end{lemma}

\subsection{Elementary Graphs: Paths and Cycles}\label{elementary}

We begin the program of determining when the property described in Problem~\ref{prob:Cameron} is true by starting with the simplest graphs: paths and cycles.  We should note that though Theorem~\ref{outerplanar} takes care of odd cycles, the techniques and observations used to prove the property described in Problem~\ref{prob:Cameron} holds for paths and cycles will play a key role in proving Theorem~\ref{outerplanar}.

\begin{proposition}\label{prop:cycle}
If $\Gamma$ is a cycle or a path and $G$ is any group of automorphisms of $\Gamma$, then the real roots of $OP_{\Gamma,G}(x)$ are bounded above by the largest real root of $P_{\Gamma}(x)$.
\end{proposition}
\proof  Throughout, we use the notation $\mathcal{C}_n$ for a cycle on $n$ vertices, and $\mathcal{P}_n$ for a path on $n$ vertices.  We begin our investigation with $\mathcal{P}_n$.  The automorphism group of $\mathcal{P}_n$ is $\ZZ{2}$, so $G$ is either the trivial group or $\ZZ{2}$.  If $G$ is the trivial group, then $OP_{\mathcal{P}_n,G}(x) = P_{\mathcal{P}_n}(x)$ so we only need to consider when $G = \ZZ{2}$.  We write $G = \{e,g\}$ where $e$ is the identity element and $g$ is the element of order $2$.  We know $\mathcal{P}_n/e = \mathcal{P}_n$, and

\[
\mathcal{P}_n / g = \begin{cases} \mathcal{P}_{\frac{n}{2}} \mbox{ with a loop at one end} & \mbox{if } n \mbox{ is even} \\ \mathcal{P}_{\frac{n - 1}{2}} & \mbox{if } n \mbox{ is odd} \end{cases}
\]
When $n$ is even, $\mathcal{P}_n/g$ has a loop so $OP_{\mathcal{P}_n,G}(x) = P_{\mathcal{P}_n}(x)$ and we are done.  If $n$ is odd, the real roots of $P_{\mathcal{P}_{\frac{n-1}{2}}}(x)$ are $\{0,1\}$, and since the same is true of $P_{\mathcal{P}_n}(x)$, the result follows by Corollary~\ref{cor:reduction}.  This successfully establishes the veracity of the property described in Problem~\ref{prob:Cameron} for paths.

We move on to cycles $\mathcal{C}_n$ where we can assume $n \geq 3$.  Now $P_{\mathcal{C}_n}(x) = (x-1)^n + (-1)^n(x-1)$ so the roots of $P_{\mathcal{C}_n}(x)$ are $\{0,1\}$ if $n$ is even, and $\{0,1,2\}$ otherwise.  We exploit this throughout our investigations with the graphs $\mathcal{C}_n$.

Let $\mathcal{C}_n$ have vertex set $\{1,2,3,\ldots,n\}$.  The automorphism group of $\mathcal{C}_n$ is the $2n$-element dihedral group, which we will denote by $D_{2n}$, whose elements are
\[
1,r,r^2,\ldots,r^{n-1},f,fr,fr^2,\ldots,fr^{n-1}
\]
where $r$ is the group element that maps $v \to v+1$ (taken mod $n$) for every $v \in V(\mathcal{C}_n)$, and $f$ is any element of order $2$ that fixes at least one vertex.  A quick computation shows

 \[\mathcal{C}_{n} / g = \begin{cases} \mbox{a loop or loops on a single vertex} & \mbox{if } g = r^{i} \mbox{ with } \gcd (i, n) = 1 \\ \mathcal{C}_{\gcd (i, n)} & \mbox{if } g = r^{i} \mbox{ with } \gcd (i, n) \neq 1 \\ \mathcal{P}_{\frac{n}{2} + 1} & \mbox{if } g = f r^{i} \mbox{ with } n \mbox{ even and } i \mbox { even} \\ \mathcal{P}_{\frac{n}{2}} \mbox{ with loops at both ends} & \mbox{if } g = f r^{i} \mbox{ with } n \mbox{ even and } i \mbox{ odd} \\ \mathcal{P}_{\frac{n + 1}{2}} \mbox{ with a loop at one end} & \mbox{if } g = f r^{i} \mbox{ with } n \mbox{ odd }\end{cases}\]

If $n$ is odd, the real roots of $P_{\mathcal{C}_n}(g)$ are $\{0,1,2\}$, and since the real chromatic roots of $\mathcal{C}_n/g$ are bounded above by $2$, the result follows by Corollary~\ref{cor:reduction}.  For the remainder of our proof, suppose $n$ is even.  

We consider all potential subgroups $G$ of $D_{2n}$.  The case when $G$ is the trivial group is immediate.  Suppose $G$ is generated by some rotation $r^i$, and without loss of generality that $i$ divides $n$.  Any non-identity element $g \in G$ will be of the form $r^{ij}$ for some positive integer $j$.  If $i$ is even, then $2|\gcd(ij,n)$ for any positive integer $j$, so each graph $\mathcal{C}_n/g$ will be an even cycle, and hence have real chromatic roots $\{0,1\}$.  Applying Corollary~\ref{cor:reduction}, the result follows.  Now if $i$ is odd, we partition the group generated by $r^i$ into sets $H_k = \{r^{(2k+1)i},r^{2ki}\}$ (exponents taken mod $n$).  If $\gcd((2k+1)i,n)=1$, then $\mathcal{C}_n/r^{(2k+1)i}$ is either a loop or a single vertex so its only real chromatic root is $0$.  The graph $\mathcal{C}/r^{2ki}$ is an even cycle so its real chromatic roots are $\{0,1\}$. By a similar argument as in the proof of Lemma~\ref{lem:reduction}, $\sum_{g \in H_k} P_{\mathcal{C}_n/g}(x)$ can not have real chromatic roots exceeding $1$.  It remains to consider when $\gcd((2k+1)i,n) \neq 1$.  For simplicity let $t_k = \gcd((2k+1)i,n)$ and $r_k = \gcd(2ki,n)$.  Observing that $\mathcal{C}_n/r^{(2k+1)i} = \mathcal{C}_{t_k}$, $\mathcal{C}_n/r^{2ki} = \mathcal{C}_{r_k}$  and $t_k$ and $r_k$ are odd and even respectively, we have
\[
P_{\mathcal{C}_n/r^{(2k+1)i}}(x) = (x-1)^{t_k} - (x-1), \ \ P_{\mathcal{C}_n/r^{2ki}}(x) = (x-1)^{r_k} + (x-1),
\]
and so 
\[
\sum_{g \in H_k} P_{\mathcal{C}_n/g}(x) = (x-1)^{t_k} + (x-1)^{r_k}.
\]
We then see that $\sum_{g \in H_k} P_{\mathcal{C}_n/g}(x)>0$ for $x>1$, and so $\sum_{g \in H_k} P_{\mathcal{C}_n/g}(x)$ does not have a real root exceeding one.  We conclude that for every $k$, the real roots of $\sum_{g \in H_k} P_{\mathcal{C}_n/g}(x)$ cannot exceed $1$, and hence by Lemma~\ref{lem:reduction}, the result follows.

The subgroups that remain are $\langle f \rangle$ and $\langle f,r^i \rangle$ (again without loss of generality, $i$ divides $n$).  Let $S = \{g \in G:g=fr^i \mbox{ for some } i\}$ and consider the partition $\mathcal{G}$ of $G$ consisting of the set $G \backslash S$ together with one-element sets, each containing a unique element from $S$.  For any $g \in S$, $\mathcal{C}_n/g$ is a path or has a loop, so the real roots of $P_{\mathcal{C}_n/g}(x)$ are $0$ or $1$.  The sum
\[
\sum_{g \in G \backslash S} P_{\mathcal{C}_n/g}(x)
\]
has roots bounded by $1$ as well, as we proved earlier.  The result then follows again by Lemma~\ref{lem:reduction}. \qed

%

%
%

\subsection{Outerplanar Graphs}

In order to establish Theorem~\ref{outerplanar}, we need to know the real chromatic roots of outerplanar graphs.  This is the content of the next proposition.

\begin{proposition}\label{prop:outerplanar}
Let $\Gamma$ be an outerplanar graph.  Then the real chromatic roots of $\Gamma$ are $\{0,1,2\}$ if $\Gamma$ contains an odd cycle, and $\{0,1\}$ otherwise.
\end{proposition}
\proof We can assume $\Gamma$ is connected, because the chromatic polynomial of a union of graphs is their product.  If $\Gamma$ is a tree, then its chromatic roots are $\{0,1\}$, so assume $\Gamma$ has a cycle.  Every such outerplanar graph $\Gamma$ can be constructed from a sequence of subgraphs \[\Gamma_1 \subset \Gamma_2 \subset \cdots \subset \Gamma_n = \Gamma\] where $\Gamma_1$ is a cycle, and $\Gamma_{i+1}$ is obtained from the subgraph $\Gamma_i$ by either:
\begin{enumerate}
\item adding an ear from $u \in V(\Gamma_i)$ to $v \in V(\Gamma_i)$ where $u,v$ are adjacent in $\Gamma_i$, or
\item adding a tree to $\Gamma_i$ with one vertex in common with $\Gamma_i$.
\end{enumerate}
In Case (1), Lemma~\ref{lem:cyclereduction} shows that the real chromatic roots of $\Gamma_{i+1}$ are those of $\Gamma_i$ with the potential addition of the real root $2$ if and only if the new cycle formed by adding the ear is an odd cycle.  In Case (2), $P_{\Gamma_{i+1}}(x) = (x-1)^{t} P_{\Gamma_{i}}(x)$ where $t$ is one fewer than the number of vertices in the tree being added, so the real chromatic roots of $\Gamma_{i+1}$ are those of $\Gamma_i$ with the potential addition of the root $1$.  The result then follows inductively. \qed

We can now establish Theorem~\ref{outerplanar}.

\proof (of Theorem~\ref{outerplanar})
We claim it suffices to show that if $\Gamma$ is outerplanar, then for any group of automorphisms $G$ of $\Gamma$ and any $g \in G$, the graph $\Gamma/g$ is either outerplanar or has a loop.  To see why, note that because $\Gamma$ contains an odd cycle, Proposition~\ref{prop:outerplanar} shows that its maximum real chromatic root is $2$.  However, if for any $g \in G$ we have that $\Gamma/g$ is outerplanar, then its real chromatic roots are bounded above by $2$ as well.  For all other $g \in G$, $\Gamma/g$ has a loop, so by Corollary~\ref{cor:reduction}, it follows that the real roots of $OP_{\Gamma,G}(x)$ are bounded above by the largest real root of $P_{\Gamma}(x)$.

It therefore remains to show that for any $g \in G$, $\Gamma/g$ is outerplanar or has a loop, provided $\Gamma$ is outerplanar itself.  We first prove this when $\Gamma$ is $2$-connected.  In this case, $\Gamma$ has a Hamiltonian cycle that forms the unique outer face of $\Gamma$ (see \cite{HararyChartrand}).  Suppose this cycle has vertices $\{1,2,3,\ldots,n\}$ in that order (and hence $V(G)=\{1,2,\ldots,n\}$).  Since the Hamiltonian cycle is unique, it must map to itself, so the group $G$ must be a subgroup of $D_{2n}$, so $g=r^i$ or $g = fr^i$ where $i \in \{0,1,2,\ldots,n-1\}$ (see Section~\ref{elementary} for definitions).  

If $g=fr^i$ for some $i$, then $\Gamma/g$ does not contain a loop only if $n$ is even and $g$ is the flip across the axis through some vertex $j$ and $\frac{n}{2}+j$ (taken mod $n$).  Consider the subsets $V_1 = \{j,j+1,\ldots,j+\frac{n}{2}\}$, $V_2 = \{j+\frac{n}{2},j+\frac{n}{2}+1,\ldots,j\}$ (taken mod $n$).  There can not be an edge $v_1v_2$ in $\Gamma$ with $v_i \in V_i$ (except possibly an edge whose endpoints are $j$ and $j+\frac{n}{2}$), for otherwise, the edge $g(v_1)g(v_2)$ (which is necessarily in $\Gamma$) would cross $v_1v_2$, contradicting the planarity of $\Gamma$.  Thus, $\Gamma/g$ is the induced subgraph of $\Gamma$ on $V_1$, with the addition of an edge from $j$ to $j+\frac{n}{2}$.  This is outerplanar since $\Gamma$ is.  

Now suppose $g=r^i$.  The orbits of $r^i$ are the same as that of $r^{\gcd(n,i)}$, so $\Gamma/r^i=\Gamma/r^{\gcd(n,i)}$, and so we work with $\Gamma/r^{\gcd(n,i)}$.  For simplicity let $k=\gcd(n,i)$.  Among all longest chords in $\Gamma$, pick the one $j_1j_2$, where $j_1<j_2$ and $j_1$ is minimal.  Observe $j_2-j_1 \leq k$, for otherwise, taking the full orbit of this chord under the group generated by $r^k$ will result in two intersecting chords in $\Gamma$, contradicting the planarity of $\Gamma$.  We then have that $\Gamma/r^k$ is the induced graph of $\Gamma$ on the vertices $\{j_1,j_1+1,\ldots,j_1+k\}$, with the potential addition of the edge from $j_1$ to $j_1+k$, and the potential addition of loops.  This is outerplanar if it doesn't have loops, and negligible if it has loops.

%
Finally, if $\Gamma$ is not $2$-connected, its biconnected components must map to each other under the action of $g$, and bridges must map to bridges, so $\Gamma/g$ identifies isomorphic biconnected components (or bridges), and remains outerplanar or has a loop.

\qed


\section{Open Problems}

It still remains to determine when the property described in Problem~\ref{prob:Cameron} holds in general.  From Example~\ref{ex2} we see that this it does not hold for all planar graphs, however Theorem~\ref{outerplanar} establishes a large class of planar graphs for which the property does hold. This leaves us with the following problem.

\begin{problem}
Characterize the planar graphs $\Gamma$ and groups $G$ for which the real roots of $OP_{\Gamma,G}(x)$ are bounded above by the largest real root of $P_{\Gamma}(x)$.
\end{problem}

Another point of interest is comparing the spread between orbital chromatic and chromatic roots.  Though we know that the real roots of $OP_{\Gamma,G}(x)$ can be larger than those of $P_{\Gamma}(x)$, how far apart can these roots be?  Based on limited experimentation in this light, we conjecture the following:

\begin{conjecture}\label{conjecture2}
For any $N>0$, there exists a graph $\Gamma$ and automorphism group $G$ of $\Gamma$ for which $OP_{\Gamma,G}(x)$ has a root at least $N$ larger than the largest real root of $P_{\Gamma}(x)$.
\end{conjecture}

The proof of Theorem~\ref{counterexample} suggests a potentially viable approach to proving Conjecture~\ref{conjecture2}: it is sufficient to find a graphs $\Gamma$ with automorphism groups $G$ for which $|V(\Gamma/g)|<|V(\Gamma/h)|$ for any $h \in G, h \neq g$, and where the largest real root of $P_{\Gamma/g}(x)$ is arbitrarily larger than that of $P_{\Gamma}(x)$.  

\section{Acknowledgments}

The authors express sincere thanks to the anonymous referees for providing a number of helpful suggestions to improve the presentation of this paper.

\bibliographystyle{plain}
\bibliography{references}

\end{document}